\documentclass[10pt,reqno]{amsart}
\usepackage{amssymb,mathrsfs,color}
\allowdisplaybreaks
\usepackage{enumerate}
\usepackage{graphicx}
\usepackage{xcolor} 
\usepackage{geometry}
\usepackage{amsfonts,amssymb,amsmath,mathtools}
\usepackage{slashed}
\usepackage{ mathrsfs }
\usepackage[titletoc,title]{appendix}
\usepackage{wrapfig}

\usepackage{cancel}

\usepackage{cite}

\definecolor{green}{rgb}{0,0.8,0} 

\title[Viscous Shock Waves of the Compressible Navier-Stokes System]{Stability of Viscous Shock Waves of the Compressible Navier-Stokes System with Large Initial Perturbation}

\author[Zhang]{Lan Zhang}
\address[Lan Zhang]{\newline School of Mathematics and Statistics and Center for Mathematical Sciences,
Wuhan University of Technology, Wuhan 430070, China\\
\newline Hubei Key Laboratory of Mathematics \& Intelligent Computation and Interdisciplinary Sciences, Wuhan University of Technology, Wuhan 430070, China}
\email{zhanglan91@whut.edu.cn}

\author[Zhao]{Huijiang Zhao}
\address[Huijiang Zhao]{\newline School of Mathematics and Statistics, Wuhan University, Wuhan 430072, China\\
 \newline Computational Science Hubei Key Laboratory, Wuhan University, Wuhan 430072, China}
\email{hhjjzhao@hotmail.com}

\newcommand{\bbR}{\mathbb R}

\definecolor{deepgreen}{cmyk}{1,0,1,0.5}

\makeatletter

\newcommand{\Rmnum}[1]{\expandafter\@slowromancap\romannumeral #1@}
\makeatother

\newcommand{\Del}[1]{}

\numberwithin{equation}{section}

\newtheorem{theorem}{Theorem}[section]
\newtheorem{lemma}[theorem]{Lemma}
\newtheorem{proposition}[theorem]{Proposition}
\newtheorem{remark}[theorem]{Remark}

\begin{document}
\date{\today}

\thanks{
\textbf{Acknowledgment.} The research of Lan Zhang was supported by the grants from National Natural Science Foundation of China (Grants No. 12101472 and 12471210).
The research of Huijiang Zhao is supported by the National Natural Science Foundation of China (Grants No. 12221001, 12371225, and 12571242).}


\begin{abstract}
We show in this paper that viscous shock waves of the one-dimensional compressible Navier-Stokes equations are time-asymptotically nonlinear stable with large initial perturbation. While results for small initial data are well-established, the case of large initial perturbations remains less satisfactory. As far as we know,
only under a class of large initial perturbations, the asymptotical stability of weak viscous shock wave was derived.
In this work, we will provide a complete resolution of this problem for weak viscous shock wave.
The key idea involves combining the $a$-contraction with shifts theory with a time splitting method to obtain the weighted $L^2$ estimates, which enables us to derive uniform upper and lower bounds on the specific volume. Our method can be also applied to proving the nonlinear stability of composite waves.
\\[2mm]
\\[2mm]
\end{abstract}

\maketitle \centerline{\date}



\section{Introduction and main results}
\subsection{The models}

In this paper, we will consider the one-dimensional isentropic compressible Navier-Stokes equations in Lagrangian coordinates as follows:
\begin{equation}\label{ns}
  \begin{cases}
  v_t-u_{x}=0,\\
  u_t+p(v)_x=\left(\frac{\mu u_x}{v}\right)_x.
  \end{cases}
\end{equation}
Here, $t\geq0$ and $x\in\bbR$. The unknown functions $v(t,x)>0$ and $u(t,x)$ stand for the specific volume and velocity, respectively. We take $p(v)=bv^{-\gamma}$ with constants $\gamma>1$ and $b>0$ which can be normalized to be $1$ without loss of generality. $\mu>0$ is the viscosity coefficient.

The main purpose of this paper is concerned with the construction of global strong solutions and the precise description of the large time behaviors of the global solutions constructed to the following Cauchy problems with prescribed initial data
\begin{equation}\label{initial-data}
  \left(v(0,x), u(0,x)\right)=\left(v_0(x), u_0(x)\right),
\end{equation}
where the initial data $\left(v(x), u(x)\right)$ is assumed to satisfy the following far-field conditions
\begin{eqnarray}\label{boundary value}
  \lim_{x\rightarrow\pm\infty}\left(v_{0}(x), u_{0}(x)\right)&=&\left(v_{\pm},u_{\pm}\right).
\end{eqnarray}
Here $v_\pm>0$ and $u_\pm,$ are constants satisfying $(v_{+},u_{+})\neq(v_{-},u_{-})$.

The large time behaviors of the global solution $(v(t,x), u(t,x))$ of the Cauchy problem \eqref{ns}, \eqref{initial-data}, \eqref{boundary value} can be determined by the structure of the unique global entropy solution of the resulting Riemann problem of the following Euler system
\begin{eqnarray}\label{approx-o-f}
  v_t-u_x&=&0,\nonumber\\
  u_t+p(v)_x&=&0,
\end{eqnarray}
with Riemann initial data
\begin{equation}\label{Riemann-data}
  (v(0,x),u(0,x))=
  \begin{cases}
  (v_-,u_-),\quad x<0,\\
  (v_+,u_+),\quad x>0.
  \end{cases}
\end{equation}

The well-known solution of the Riemann problem for the system \eqref{approx-o-f} and \eqref{Riemann-data}, first proposed and
solved by Riemann, which is determined by a combination of at most two elementary solutions from the following
four families: 1-rarefaction; 2-rarefaction; 1-shock and 2-shock. In this paper, we only focus on the shock waves.
The shock curves are defined by using Rankine-Hugoniot condition,
there exists a constant $s$ such that
\begin{align*}
  &-s(v_+-v)-(u_+-u)=0,\\
  &-s(u_+-u)+(p(v_+)-p(v))=0,
\end{align*}
which yields
\begin{equation*}
  s=\sqrt{-\frac{p(v_+)-p(v_-)}{v_+-v_-}},
\end{equation*}
for the the 2-shock wave $S_2(v_+,u_+)$ with admissible shocks for $0<v_-<v< v_+$. The solution to system \eqref{approx-o-f} and \eqref{Riemann-data}
is given by the discontinuous traveling wave defined as
\begin{equation*}
  (v^s(t,x),u^s(t,x))=
  \begin{cases}
  (v_-,u_-),\quad x<st,\\
  (v_+,u_+),\quad x>st.
  \end{cases}
\end{equation*}

Next, we introduce the viscous shock waves. By introducing a new variable $\xi=x-st$,  the the 2-viscous shock wave $(V,U)(\xi)$ satisfies the following ODE:
\begin{equation}\label{shock-wave}
  \begin{cases}
  -sV_{\xi}-U_{\xi}=0,\\
  -sU_{\xi}+p(V)_{\xi}=\left(\frac{\mu U_{\xi}}{V}\right)_{\xi},\\
  \lim_{\xi\pm\infty}(V,U)=(v_{\pm},u_{\pm}).
  \end{cases}
\end{equation}

\subsection{The main results}
Under the above preparations, our first result on the Cauchy problem \eqref{ns}, \eqref{initial-data}, \eqref{boundary value} for the Navier-Stokes equations can be stated as in the following theorem.
\begin{theorem}\label{theorem1}
Let $\delta:=|v_+-v_-|$ be the strength of the shock wave solution $\left(V,U\right)(x-st-X(t))$ of the second family of the Riemann problem \eqref{approx-o-f}, \eqref{Riemann-data} and assume that $(v_-,u_-)\in S_2(v_+,u_+)$. If we assume further that
the initial data $(v_0(x), u_0(x))$ satisfy
\begin{eqnarray}\label{initial-norm1}
  \left(v_{0}(x)-V_0(x), {u}_{0}(x)-U_0(x)\right)&\in& L^2(\bbR),\nonumber\\
  \left(\partial_x v_0(x), \partial_x u_0(x)\right)&\in& L^2(\bbR),\\
  \inf_{x\in\bbR}\big\{v_{0}(x)\big\}&>&0,\nonumber
\end{eqnarray}
then there exists a sufficiently small positive constant $\delta_0$ such that if $0<\delta\leq \delta_0$, the Cauchy problem
 \eqref{ns}, \eqref{initial-data}, \eqref{boundary value} for the Navier-Stokes equations admits a unique global strong
solution $(v(t,x), u(t,x))$, which satisfies
\begin{equation}\label{bounds-1}
C_1^{-1}\leq v(t,x)\leq C_1,
\end{equation}
for all $(t,x)\in\bbR_+\times \bbR$. Moreover, there exists an absolutely continuous shift $X(t)$ such that
\begin{equation}\label{thm1-asmptotic}
  \lim_{t\rightarrow+\infty}\sup_{x\in\bbR}\left|
  \left(v(t,x)-V(x-st-X(t)),u(t,x)-U(x-st-X(t))\right)\right|=0.
\end{equation}
and
\begin{equation}\label{X-asmp}
  \lim_{t\rightarrow+\infty}|\dot{X}(t)|=0.
\end{equation}
Here $\delta_0$, $C_1\geq1$, and $C_2$ are some positive constants depending only on $\inf\limits_{x\in\bbR}\{n_{0}(x)\}$, $\|(\partial_x v_0, \partial_x u_0)\|_{L^2(\mathbb{R})}$, and
$\left\|\left({v}_{0}-V_0, {u}_{0}-U_0\right)\right\|_{L^2(\mathbb{R})}$.
\end{theorem}

\begin{remark}
This is the first result concerning the asymptotical stability of viscous shock waves for one dimensional isentropic Navier-Stokes equations with
general large initial data.
The method here can be modified directly to the composite waves of 1-rarefaction wave and 2-viscous shock wave with shift $X(t)$. For one
dimensional non-isentropic Navier-Stokes equations, we can also obtain that the composite waves of 1-rarefaction wave, 2-contact
discontinuity and 3-shock wave are time-asymptotically nonlinear stable with large initial data. One can find these results in our forthcoming papers.
\end{remark}
\begin{remark}
Recently, there are excited results concerning the decay rates of these fundamental waves when the initial data are sufficiently small,
see \cite{CW-2026,Chen-2026,HX-2026,HHLX-2026}.
In a companion paper, we will show that such decay rates are also established for large initial data.
\end{remark}

\subsection{Former related results and main ideas}
It is well-known that the nonlinear stability of shock waves for one dimensional Navier-Stokes equations with small initial data are well-established.
Matsumura-Nishihara \cite{MR839317} proved such results by using the so called anti-derivative method under the zero mass condition. Same methods are
independently applied by Goodman \cite{MR853782} for general conservation laws with artificial diffusion. The crucial zero mass condition was then removed
by Liu \cite{MR791863}, Szepessy-Xin \cite{MR1207241} and Liu-Zeng \cite{MR3244333}. A spectral stability of viscous shock to 1D compressible Navier-Stokes system under
a spectral condition was proved by \cite{MR2044067}. Recently, a favorable method called ``$a$-contraction with shifts'' was developed
by Kang-Vasseur \cite{MR4195742,MR4228501}. Depending on this method, a conjecture mentioned by Matsumura-Nishihara \cite{MR899210} concerning the time-asymptotic stability of composite waves of viscous shock waves and rarefaction wave for barotropic Naiver-Stokes equations was accomplished in \cite{MR4560999}. The
corresponding results for compressible Navier-Stokes-Fourier system were proved in \cite{MR4924422}.

However, the similar results for large initial data are less satisfactory. For isentropic case, Wang-Zhao-Zou \cite{MR3063535} deduce the nonlinear stability result
for certain class of large initial perturbation which can allow the initial density to have large oscillation. More precisely, the initial perturbations satisfy
\begin{equation*}
  \|(\phi_0,\psi_0)\|_{L^2}\leq \delta^{\alpha},\quad \|(\phi_{0x},\psi_{0x})\|_{L^2}\leq \delta^{-\beta},
\end{equation*}
where $\beta>\alpha>0$, hence the norm $\|(\phi_0,\psi_0)\|_{L^{\infty}}$ can be arbitrarily large. Similar results are extended to non-isentropic case in \cite{MR3861291,MR3794193}. However, all the above results need to assume that the $L^2$ norm of the initial perturbations are sufficiently small. This crucial
assumption is used to close the $L^2$ type energy estimates, which enables them to derive uniformly upper and lower bounds on the specific volume $v(t,x)$.

In this paper, we decide to remove such assumption to obtain the time-asymptotically nonlinear stability of weak shock waves.
Based on above analysis, we conclude that the main difficult lies in closing the $L^2$ type energy estimates.
Recently, Kang-Vasseur \cite{MR4195742,MR4228501} studied a contraction property for viscous weak shocks, up to a shift, for any possibly large perturbations.
In their article, a weighted $L^1$ type estimate for the relative entropy is developed without any smallness condition imposed on the initial data,
which is independent of the strength of the viscosity.
These results motivate us to consider the possibility to extend the contraction property to large time behavior of viscous shock waves.
However, the estimates in \cite{MR4195742,MR4228501} can not be applied directly to our case, due to the crucial assumptions on the viscosity that is $\mu=bv^{-\alpha}$,
with $0<\alpha\leq \gamma\leq \alpha+1$. Recently, \cite{EEKE-2025} improved the results for any $\alpha\in[0,1]$ for compressible isothermal Navier-Stokes equations,
that is $p(v)=\frac{1}{v}$. This result is still not suitable for our case since we consider $p(v)=v^{-\gamma}$ with $\gamma>1$. Hence, we need some new ideas.

The main idea is splitting the temporal estimates of nonlinear terms into following two parts:
\begin{equation*}
  \mathfrak{G}(t):=\{t\in[0,T]:\,\|(\phi,\psi)(t,x)\|_{L^{\infty}}\leq \eta\},\qquad \mathfrak{B}(t):=\{t\in[0,T]:\,\|(\phi,\psi)(t,x)\|_{L^{\infty}}\geq \eta\}.
\end{equation*}
In the good time set $\mathfrak{G}$, we can follow the ideas constructed in \cite{MR4195742,MR4228501,MR4560999} to control the possible growth.
In the bad time set $\mathfrak{B}$, although the point-wise bounds of $(\phi,\psi)$ can be large, we can connect it with the $L^2$ estimates of its gradient
as follows:
\begin{align*}
  1_{\mathfrak{B}(t)}\leq \eta^{-4}\|(\phi,\psi)\|_{L^{\infty}}^4\leq C(M_1)\eta^{-4}N_2\|(\phi_x,\psi_x)\|_{L^2}^2.
\end{align*}
Then, we find that all the nonlinear terms are bounded by
\begin{equation*}
  (C(M_1)(1+k^2)+c_2)\int_0^T1_{\mathfrak{B}}\int_{\bbR}V^X_{\xi}(\phi^2+\psi^2)dxdt.
\end{equation*}
Since
\begin{equation*}
  \int_{\bbR}V^X_{\xi}(\phi^2+\psi^2)dx\leq C(M_1)\delta^2 N_2,
\end{equation*}
combining the above estimates with \eqref{first-order-est}, we finally close the weighted $L^2$ type estimates after fixing the parameters $M_1$, $k$, $\eta$ and $N_2$
and choosing the strength of the shock waves $\delta$ sufficiently small. With this estimates, we can then derive the upper and lower bounds for the specific volume and
the $H^1$ estimates. Applying local existence theorem and bootstrap analysis, the Theorem \ref{theorem1} is thus proved.

\subsection{Outline of this paper and notations}
The rest of the paper is organized as follows: We prepare some inequalities and useful constructions for proof in Section 2.
Section 3 is devoted to deriving the energy estimates and uniformly upper and lower bounds on the specific volume.
The Theorem \ref{theorem1} is finally proved in Section 4.\\

\noindent \textit{Notations.} For notational simplicity, throughout this paper, $C$ denotes some generic positive constant independent of the parameter $0<\delta\leq 1$, which may take different values in different places. $C(\varepsilon)$ and $C(M_i)$ represent some positive constants depending only on $\varepsilon$ and $M_i$, respectively. $A\lesssim B$ means that there is a generic constant $C>0$ independent of $0<\delta\leq 1$ such that $A\leq CB$. $A\thicksim B$ means $A\lesssim B$ and $B\lesssim A$.
$L^q(\mathbb{R})(1\leq q\leq \infty)$ stands for the usual Lebesgue space with norm $\|\cdot\|_{L^q}$,
and $H^k(\mathbb{R})(k\in \mathbb{Z}_+)$ is the usual Sobolev space with norm $\|\cdot\|_{H^k}$.

\section{Preliminary}
\subsection{Properties of viscous shock wave}
The properties of the 2-viscous shock wave $(V,U)(\xi)$ can be listed in the following lemma.
\begin{lemma}[cf. \cite{MR839317}]\label{shock-property}
 For any state $(v_+, u_+)$, there exists a constant $C > 0$ such that the following
is true. For any end state such that $(v_-, u_-)\in S_2(v_+, u_+)$, there exists a unique solution $(V,U)(\xi)$ to \eqref{shock-wave}.
We fix the normalization $V(0)=\frac{v_++v_-}{2}$.
 Let $\delta$ denote the strength of the shock as $\delta=v_+-v_-$. It holds that
 \begin{equation*}
   V_{\xi}>0,\quad  U_{\xi}=-sV_{\xi}<0.
 \end{equation*}
 and
 \begin{align*}
   &|(V_{\xi},U_{\xi})|\leq C\delta^2e^{-C\delta|\xi|},\\
   &|(V_{\xi\xi},U_{\xi\xi})|\leq C\delta|(V_{\xi},U_{\xi})|.
 \end{align*}
\end{lemma}

\subsection{Useful inequality}
The method of a-contraction with shift relies on the following Poincar\'{e} type inequality
\begin{lemma}[cf. \cite{MR4195742}]\label{Poincare-ineq}
For any $f:[0,1]\rightarrow\bbR$, satisfying $\int_0^1y(1-y)|f'|^2dy<\infty$,
\begin{equation*}
  \int_0^1|f-\int_0^1fdy|^2dy\leq \frac12\int_0^1y(1-y)|f'|^2dy.
\end{equation*}
\end{lemma}

\subsection{Construction of weight function}
We define a weight function $a^X(t,x)$ by
\begin{equation}\label{a}
   a^X=1+\frac{\lambda}{\delta}(V^X-v_-),
\end{equation}
where we define $k:=\frac{\lambda}{\delta}$ which satisfies $k\gg 1$ such that the constant $\lambda=k\delta\ll 1$.

Notice that
\begin{equation*}
  a^X_{\xi}=kV^X_{\xi}>0, \quad a^X\in[1,2].
\end{equation*}


\subsection{Construction of shift}
We define the shift $X$ as a solution to the ODE:
\begin{equation}\label{shift}
  \begin{cases}
  \dot{X}(t)=-\frac{p''(v_-)}{(-p'(v_-))^{\frac32}}\frac{1}{\delta}\int_{\bbR}a^X(U^X_{\xi}\psi-p'(V^X)V^X_{\xi}\phi)dx,\\
  X(0)=0
  \end{cases}
\end{equation}
The system \eqref{shift} has a unique absolutely continuous solution defined on any interval in time $[0, T]$, whose
proof can be found in \cite[Lemma 3.3]{MR4560999}.

\section{Local existence and a priori estimates}
We consider the Cauchy problem of Navier-Stokes equations \eqref{ns} in this section and our main purpose is to prove Theorem \ref{theorem1}.
\subsection{Local existence}
Define the following perturbation
\begin{equation*}
  \left(\phi(t,x), \psi(t,x)\right)=\left(v(t,x)-V(x-st-X(t)), u(t,x)-U(x-st-X(t))\right),
\end{equation*}
where we denote $V^X(t,x)=V(x-st-X(t))$ and $U^X(t,x)=U(x-st-X(t))$.
Then we can deduce from \eqref{ns} and \eqref{shock-wave} that $\left(\phi(t,x), \psi(t,x)\right)$ solves
\begin{eqnarray}\label{perturbation-o-f}
  \phi_t-\psi_x&=&\dot{X}(t)V^X_x,\nonumber\\
  \psi_t+\left(p(v)-p(V^X)\right)_x&=&\left(\frac{\mu u_x}{v}-\frac{\mu U_x^X}{V^X}\right)_x+\dot{X}(t)U^X_x
\end{eqnarray}
with prescribed initial data
\begin{equation}\label{perturb-initial-data-0-f}
  \left(\phi(0,x),\psi(0,x)\right)=\left(\phi_0(x),\psi_0(x)\right):=\left(v_{0}(x)-V^X_0(x),u_0(x)-U^X_0(x)\right).
\end{equation}

To prove Theorem \ref{theorem1}, for some given positive constants $T, M_1\geq1, M_2$, and $N_1$, we first define the following set of functions for which the solutions to the Cauchy problem \eqref{perturbation-o-f},  \eqref{perturb-initial-data-0-f} will be sought:
\begin{eqnarray}\label{solution-space-o-f}
  &&X\left(0,T;M_1,N_1\right)\\
  &:=&\left\{\left(\phi(t,x), \psi(t,x)\right)\ \left|
\begin{array}{lc}
&\left(\phi(t,x), \psi(t,x)\right)\in C\left(0,T;H^1(\mathbb{R})\right),\\
&\partial_x\phi(t,x)\in L^2\left(0,T;L^2(\mathbb{R}\right),\\
&\partial_x\psi(t,x)\in L^2\left(0,T;H^1(\mathbb{R})\right),\\
&M_1^{-1}\leq \phi(t,x)+V^X(t,x)\leq M_1,\quad \forall (t,x)\in[0,T]\times\mathbb{R},\\
&\left\|\left(\phi,\psi\right)\right\|_{L^{\infty}\left(0,T;H^1(\mathbb{R})\right)}\leq N_1
\end{array}
\right.
\right\}.\nonumber
\end{eqnarray}

For the local solvability of the Cauchy problem \eqref{perturbation-o-f}, \eqref{perturb-initial-data-0-f}, we have by employing the standard iteration method that
\begin{proposition}\label{local-exist}
Suppose $(\phi_{0}(x), \psi_{0}(x))\in H^1(\mathbb{R})$ and $\inf\limits_{x\in\mathbb{R}} v_{0}(x)>0$, then there exists a sufficiently small positive constant $t_1$, which depends only on $N^0_1=\|(\phi_{0}, \psi_{0})\|_{H^1}$ and $M^0_1=\min\left\{\left(\inf_{x\in\mathbb{R}} n_0(x)\right)^{-1},\|n_0\|_{L^\infty(\mathbb{R})}\right\}$, such that the Cauchy problem \eqref{perturbation-o-f}$, $\eqref{perturb-initial-data-0-f} admits a unique strong solution $(\phi, \psi) \in X\left(0,t_1;2M^0_1,2N^0_1\right)$.
\end{proposition}

\subsection{Relative entropy method}
The entropy of system \eqref{ns} is $\eta(v,u):=\frac12 u^2+Q(v)$, with $Q(v)=\frac{1}{\gamma-1}v^{1-\gamma}$. Notice that $Q'(v)=-p(v)$.
Consider the relative entropy functional defined by
\begin{equation*}
  \eta((v,u)|(V^X,U^X))=\frac12 \psi^2+Q(v|V^X).
\end{equation*}
with
\begin{equation*}
  Q(v|V^X)=Q(v)-Q(V^X)+p(V^X)(v-V^X).
\end{equation*}

Suppose that the above local solution $(\phi(t,x), \psi(t,x))$ has been extended to the time step $t=T\geq t_1$ and assume that the following {\it a priori} assumption
$$
\left(\phi(t,x), \psi(t,x)\right) \in X\left(0,T;M_1, N_1\right)
$$
holds for some positive constants $M_1$ and $N_1$, then we want to deduce the desired $H^1(\mathbb{R})-$norm {\it a priori} estimates on $(\phi(t,x), \psi(t,x))$ in terms of the initial data $\left(\phi_0(x),\psi_0(x)\right)$ but independent of $M_1, N_1$, which is the main content of the following proposition.
\begin{proposition}\label{proposition-o-f}
Let $\left(\phi(t,x), \psi(t,x)\right)\in X(0,T;M_1,N_1)$ be a solution to
the Cauchy problem \eqref{perturbation-o-f}$, $\eqref{perturb-initial-data-0-f} defined on the strip $\Pi_T=[0,T]\times\mathbb{R}$, then if $\delta=|v_+-v_-|$, the strength of the shock wave $\left(V, U\right)(x-st-X(t))$ is chosen sufficiently small such that
\begin{equation}\label{delta-assumption}
  \int_{\bbR}a^X\eta((v,u)|(V^X,U^X))dx\leq N_2,\quad \delta N_2\ll1,
\end{equation}
the solution $\left(\phi(t,x), \psi(t,x)\right)$ satisfies
\begin{align}
  \sup_{0\leq \tau\leq t}\|(\phi,\psi)(\tau)\|_1^2+\int_{0}^{t} \big(\|\partial_x\phi(\tau)\|^2+\|\partial_x\psi(\tau)\|_{H^1}^2 \big) \mathrm{d}\tau
   \lesssim 1.\label{pro-o-1}
\end{align}
for all $0\leq t\leq T$. Moreover, there exist positive constants $M_2\geq1$, which depend only on $N^0_1=\|(\phi_{0}, \psi_{0})\|_{H^1}$ and $M^0_1=\min\left\{\left(\inf_{x\in\mathbb{R}} n_0(x)\right)^{-1},\|n_0\|_{L^\infty(\mathbb{R})}\right\}$, such that
\begin{equation*}
  M_2^{-1}\leq v(t,x)\leq M_2.
\end{equation*}
\end{proposition}

The Proposition \ref{proposition-o-f} will be proved by the following Lemmas.
Similarly, we define a relative pressure
\begin{equation*}
  p(v|V^X)=p(v)-p(V^X)-p'(V^X)(v-V^X).
\end{equation*}
Below we will consider the relative entropy weighted by $a^X$ defined in \eqref{a}.
\begin{lemma}\label{lem-rela-weight-entropy}
Let $a^X$ be the weight function defined by \eqref{a}. Let $(v,u)$ be a solution to \eqref{ns}. Then
\begin{equation}\label{rela-entropy-est}
  \frac{d}{dt}\int_{\bbR}a^X\eta((v,u)|(V^X,U^X))dx=\dot{X}(t)Y(v,u)-\mathcal{J}^{good}+\mathcal{J}^{bad},
\end{equation}
with
\begin{align*}
  Y(v,u)=&\,\int_{\bbR}a^X(U^X_{\xi}\psi-p'(V^X)V^X_{\xi}\phi)dx-\int_{\bbR}a^X_{\xi}\eta((v,u)|(V^X,U^X))dx,\\
  \mathcal{J}^{good}=&\,\int_{\bbR}\frac{\mu}{v} a^X\psi_x^2dx+\frac{\sqrt{-p'(v_-)}}{2}\int_{\bbR}a^X_{\xi}(\sqrt{-p'(v_-)}\phi+\psi)^2dx,\\
  \mathcal{J}^{bad}=&\,s\int_{\bbR}a^XV^X_{\xi}p(v|V^X)dx-\frac{1}{2}(s-\sqrt{-p'(v_-)})\int_{\bbR}a^X_{\xi}(\psi^2-p'(v_-)\phi^2)dx\\
  &\,+\int_{\bbR}a^X_{\xi}\left[\phi\psi\int_0^1(p'(V^X+\tau\phi)-p'(v_-))d\tau-s\phi^2\int_0^1(1-\tau)(p'(V^X+\tau\phi)-p'(v_-))d\tau\right]dx\\
  &\,-\int_{\bbR}\mu a^XU^X_{\xi}\psi_x\left(\frac{1}{v}-\frac{1}{V^X}\right)dx-\int_{\bbR}\mu a^X_{\xi}\psi\left(\frac{\psi_x}{v}+U^X_{\xi}\left(\frac{1}{v}-\frac{1}{V^X}\right)\right)dx.
\end{align*}
\end{lemma}
\begin{proof}
Direct calculations show that
\begin{align}\label{rela-entropy-identity}
  &\,\frac{d}{dt}\int_{\bbR}a^X\eta((v,u)|(V^X,U^X))dx\notag\\
  =&\,\int_{\bbR}a^X_{\xi}(-s-\dot{X}(t))\eta((v,u)|(V^X,U^X))dx+\int_{\bbR}a^X\left(\psi\psi_t+\frac{d}{dt}Q(v|V^X)\right)dx\notag\\
  =&\,-\int_{\bbR}\frac{\mu}{v} a^X\psi_x^2dx-\int_{\bbR}\mu a^XU^X_{\xi}\psi_x\left(\frac{1}{v}-\frac{1}{V^X}\right)dx-\int_{\bbR}\mu a^X_{\xi}\psi\left(\frac{\psi_x}{v}+U^X_{\xi}\left(\frac{1}{v}-\frac{1}{V^X}\right)\right)dx\notag\\
  &\,+\int_{\bbR}a^X_{\xi}\left[\psi(p(v)-p(V^X))-s\eta((v,u)|(V^X,U^X))\right]dx+s\int_{\bbR}a^XV^X_{\xi}p(v|V^X)dx\notag\\
  &\,+\dot{X}(t)\left[\int_{\bbR}a^X(U^X_{\xi}\psi-p'(V^X)V^X_{\xi}\phi)dx-\int_{\bbR}a^X_{\xi}\eta((v,u)|(V^X,U^X))dx\right].
\end{align}
Next we apply mean value theorem to get
\begin{align*}
  &\,\psi(p(v)-p(V^X))-s\eta((v,u)|(V^X,U^X))\\
  =&\,\left(p'(v_-)+\int_0^1(p'(V^X+\tau\phi)-p'(v_-))d\tau\right)\phi\psi\\
  &\,-\frac{s}{2}\psi^2+\frac{s}{2}p'(v_-)\phi^2+s\phi^2\int_0^1(1-\tau)(p'(V^X+\tau\phi)-p'(v_-))d\tau\\
  =&\,-\frac{\sqrt{-p'(v_-)}}{2}\int_{\bbR}a^X_{\xi}(\sqrt{-p'(v_-)}\phi+\psi)^2dx-\frac{1}{2}(s-\sqrt{-p'(v_-)})(\psi^2-p'(v_-)\phi^2)\\
  &\,+\phi\psi\int_0^1(p'(V^X+\tau\phi)-p'(v_-))d\tau-s\phi^2\int_0^1(1-\tau)(p'(V^X+\tau\phi)-p'(v_-))d\tau.
\end{align*}
With above identities, we derive the desired formula \eqref{rela-entropy-est}.
\end{proof}

To prove Proposition \ref{proposition-o-f}, we divide the proof into two parts with respect to the time. We define two sets by good time set and bad time set as follows:
\begin{equation*}
  \mathfrak{G}(t):=\{t\in[0,T]:\,\|(\phi,\psi)(t,x)\|_{L^{\infty}}\leq \eta\},\qquad \mathfrak{B}(t):=\{t\in[0,T]:\,\|(\phi,\psi)(t,x)\|_{L^{\infty}}\geq \eta\}.
\end{equation*}
Here, the positive constant $\eta$ will be chosen later which is small and the norm $L^{\infty}$ is considered with respect to spatial variable.
In the good set, for example taking $\eta\leq \frac{v_+}{8}$, the specific volume $v(t,x)$ is hence bounded by $[\frac{3}{8}v_+,\frac{9}{8}v_+]$.
We then have the following lemma.

\begin{lemma}\label{lem-gb-time}
Under the assumptions in Proposition \ref{proposition-o-f}, there exists a generical constant $c_0$ such that
\begin{align}\label{gb-time-est}
  &\,\int_{\bbR}a^X\eta((v,u)|(V^X,U^X))dx+c_0\int_0^T\left(\int_{\bbR}\frac{\mu}{v}a^X\psi_x^2dx+\int_{\bbR}V^X_{\xi}(\phi^2+\psi^2)dx+\delta|\dot{X}(t)|^2\right)dt\notag\\
  \leq&\, \int_{\bbR}a^X\eta((v,u)|(V^X,U^X))(0,x)dx+C(M_1)(1+k^2)\eta^{-4}\delta^2N^2_2\int_0^T\|(\phi_x,\psi_x)\|_{L^2}^2dt.
\end{align}
\end{lemma}
\begin{proof}
Integrating \eqref{rela-entropy-identity} with respect to time $t$, we have
\begin{align}\label{rela-entropy-integ}
  &\,\int_{\bbR}a^X\eta((v,u)|(V^X,U^X))dx-\int_{\bbR}a^X\eta((v,u)|(V^X,U^X))(0,x)dx\notag\\
  =&\,-\int_0^T\int_{\bbR}\frac{\mu}{v} a^X\psi_x^2dxdt-\int_0^T\int_{\bbR}\mu a^XU^X_{\xi}\psi_x\left(\frac{1}{v}-\frac{1}{V^X}\right)dxdt-\int_0^T\int_{\bbR}\mu a^X_{\xi}\psi\left(\frac{\psi_x}{v}+U^X_{\xi}\left(\frac{1}{v}-\frac{1}{V^X}\right)\right)dxdt\notag\\
  &\,+\int_0^T\int_{\bbR}a^X_{\xi}\left[\psi(p(v)-p(V^X))-s\eta((v,u)|(V^X,U^X))\right]dxdt+s\int_0^T\int_{\bbR}a^XV^X_{\xi}p(v|V^X)dxdt\notag\\
  &\,+\int_0^T\dot{X}(t)\left[\int_{\bbR}a^X(U^X_{\xi}\psi-p'(V^X)V^X_{\xi}\phi)dx-\int_{\bbR}a^X_{\xi}\eta((v,u)|(V^X,U^X))dx\right]dt.
\end{align}

To estimate the right hand side of above identity, we divide the integration interval with respect to time into two parts, that is, $\int_0^T\ast dt=\int_{\mathfrak{G}}\ast dt+\int_{\mathfrak{B}}\ast dt$. In the good time set $\mathfrak{G}$, we will apply the decomposition in \eqref{rela-entropy-est}. While, in the bad time set $\mathfrak{B}$, we will directly estimate the terms in \eqref{rela-entropy-integ}.
\medbreak

\noindent{\bf 1. Good time estimates.}

Firstly, we consider the terms in $\mathcal{J}^{bad}:=\sum_{i=1}^7\mathcal{J}^{bad}_i$. For the term $\mathcal{J}^{bad}_1$, notice that
\begin{align*}
  sp(v|V^X)=&\,\frac{1}{2}\sqrt{-p'(v_-)}p''(v_-)\phi^2+\frac12(s-\sqrt{-p'(v_-)})p''(v_-)\phi^2\\
  &\,+s\phi^2\int_0^1(1-\tau)(p''(V^X+\tau\phi)-p''(v_-))d\tau.
\end{align*}
and
\begin{align*}
  &s^2+p'(v_-)=-\int_0^1(p'(v_-+\tau\delta)-p'(v_-))d\tau\, \Rightarrow\, s-\sqrt{-p'(v_-)}\leq C\sup_{}|p''|\delta,\\
  &|p'(V^X+\tau\phi)-p'(v_-)|\leq \sup_{}|p''|(\delta+\eta),\\
  &|p''(V^X+\tau\phi)-p''(v_-)|\leq \sup_{}|p'''|(\delta+\eta).
\end{align*}
We obtain
\begin{align}\label{J-bad-1}
  \int_{\mathfrak{G}}\mathcal{J}^{bad}_1dt=&\,\frac{1}{2}\sqrt{-p'(v_-)}p''(v_-)\int_{\mathfrak{G}}\int_{\bbR}a^XV^X_{\xi}\phi^2dxdt
  +\frac12(s-\sqrt{-p'(v_-)})p''(v_-)\int_{\mathfrak{G}}\int_{\bbR}a^XV^X_{\xi}\phi^2dxdt\notag\\
  &\,+\int_{\mathfrak{G}}\int_{\bbR}a^XV^X_{\xi}s\phi^2\int_0^1(1-\tau)(p''(V^X+\tau\phi)-p''(v_-))d\tau dxdt\notag\\
  \leq&\,\frac{1}{2}\sqrt{-p'(v_-)}p''(v_-)\int_{\mathfrak{G}}\int_{\bbR}a^XV^X_{\xi}\phi^2dxdt+C(\delta+\eta)\int_{\mathfrak{G}}\int_{\bbR}V^X_{\xi}(\phi^2+\psi^2)dxdt\notag\\
  \leq&\,\frac{1}{2}\sqrt{-p'(v_-)}p''(v_-)\int_{\mathfrak{G}}\int_{\bbR}a^XV^X_{\xi}\left(\frac{\psi+\sqrt{-p'(v_-)}\phi}{\sqrt{-p'(v_-)}}-\frac{\psi}{\sqrt{-p'(v_-)}}\right)^2dxdt\notag\\
  &\,\quad+C(\delta+\eta)\int_{\mathfrak{G}}\int_{\bbR}V^X_{\xi}(\phi^2+\psi^2)dxdt\notag\\
  \leq&\, (1+k^{-\frac12})\frac{p''(v_-)}{2\sqrt{-p'(v_-)}}\int_{\mathfrak{G}}\int_{\bbR}a^XV^X_{\xi}\psi^2dxdt
  +C(\delta+\eta)\int_{\mathfrak{G}}\int_{\bbR}V^X_{\xi}(\phi^2+\psi^2)dxdt\notag\\
  &\,+(1+k^{\frac12})\frac{p''(v_-)}{2\sqrt{-p'(v_-)}}\int_{\mathfrak{G}}\int_{\bbR}a^XV^X_{\xi}(\psi+\sqrt{-p'(v_-)}\phi)^2dxdt\notag\\
  \leq&\, (1+k^{-\frac12})(1+\lambda)\frac{p''(v_-)}{2\sqrt{-p'(v_-)}}\int_{\mathfrak{G}}\int_{\bbR}V^X_{\xi}\psi^2dxdt
  +Ck^{-\frac12}\int_{\mathfrak{G}}\int_{\bbR}a^X_{\xi}(\psi+\sqrt{-p'(v_-)}\phi)^2dxdt\notag\\
  &\,+C(\delta+\eta)\int_{\mathfrak{G}}\int_{\bbR}V^X_{\xi}(\phi^2+\psi^2)dxdt.
\end{align}

For the terms $\mathcal{J}^{bad}_{i}$ $(i=2\sim4)$, applying $a^X_{\xi}=kV^X_{\xi}$, we have
\begin{align}\label{J-bad-2-4}
  \sum_{i=2}^4\left|\int_{\mathfrak{G}}\mathcal{J}^{bad}_idt\right|\leq Ck(\delta+\eta)\int_{\mathfrak{G}}\int_{\bbR}V^X_{\xi}(\phi^2+\psi^2)dxdt.
\end{align}

For the terms $\mathcal{J}^{bad}_{i}$ $(i=5\sim7)$, we have
\begin{align}\label{J-bad-5-6}
  \sum_{i=5}^6\left|\int_{\mathfrak{G}}\mathcal{J}^{bad}_idt\right|\leq&\, \frac18 \int_{\mathfrak{G}}\int_{\bbR}\frac{\mu}{v}a^X\psi_x^2dxdt+C(\delta^2+k^2\delta^2)\int_{\mathfrak{G}}\int_{\bbR}V^X_{\xi}(\phi^2+\psi^2)dxdt,\\
  \left|\int_{\mathfrak{G}}\mathcal{J}^{bad}_7dt\right|\leq&\,  Ck\delta^2\int_{\mathfrak{G}}\int_{\bbR}V^X_{\xi}(\phi^2+\psi^2)dxdt.\label{J-bad-7}
\end{align}

Next, we rewrite the first term in $\mathcal{J}^{good}$ in term of the new variable $y$ and $\omega$ as follows:
\begin{align*}
  y=\frac{V^X-v_-}{\delta}\in[0,1],\quad \omega(y):=\psi(x(y)).
\end{align*}
Notice that $\frac{dy}{dx}=\frac{V^X_{\xi}}{\delta}>0$. Applying
$s\mu\frac{V^X_{\xi}}{V^X}=-s^2(V^X-v_-)-(p(V^X)-p(v_-))$
yields
\begin{equation*}
  \left|\frac{\mu y_x}{V^X(1-y)y}-\frac{p''(v_-)}{2\sqrt{-p'(v_-)}}\delta\right|\leq C\delta^2.
\end{equation*}
The same proof can be found in \cite[Lemma 3.1]{MR4195742}.
With above preparation, we have
\begin{align}\label{J-good-1}
  \int_{\mathfrak{G}}\int_{\bbR}\frac{\mu}{v}a^X\psi_x^2dxdt=&\,\int_{\mathfrak{G}}\int_0^1\frac{\mu}{v}a^X\omega_y^2y_xdydt=\int_{\mathfrak{G}}\int_0^1a^X\frac{V^X}{v}\frac{\mu y_x}{V^X(1-y)y}y(1-y)\omega_y^2dydt\notag\\
  \geq&\,\left(1-\frac{8}{3v_+}\eta\right)\left(\frac{p''(v_-)}{2\sqrt{-p'(v_-)}}\delta-C\delta^2\right)\int_{\mathfrak{G}}\int_0^1y(1-y)\omega_y^2dydt\notag\\
  \geq&\,\frac{p''(v_-)}{2\sqrt{-p'(v_-)}}\delta\left(1-C(\delta+\eta)\right)\int_{\mathfrak{G}}\int_0^1y(1-y)\omega_y^2dydt\notag\\
  \geq&\,\frac{p''(v_-)}{\sqrt{-p'(v_-)}}\delta\left(1-C(\delta+\eta)\right)\int_{\mathfrak{G}}\left(\int_0^1\omega^2dy-\left(\int_0^1\omega dy\right)^2\right)dt.
\end{align}
In the last inequality, we have used Lemma \ref{Poincare-ineq}.

Finally, for the two terms in $Y(v,u)$, we rewrite the first term in terms of new variables:
\begin{align*}
  &\,\int_{\bbR}a^X(U^X_{\xi}\psi-p'(V^X)V^X_{\xi}\phi)dx\\
  =&\,\delta\int_0^1a^X(-s\omega-p'(V^X)\phi(x(y)))dy\\
  =&\,-2\delta\sqrt{-p'(v_-)}\int_0^1\omega dy-\delta\int_0^1\left(sa^X-2\sqrt{-p'(v_-)}-\frac{a^Xp'(V^X)}{\sqrt{-p'(v_-)}}\right)\omega dy\\
  &\,-\delta\int_0^1a^X\frac{p'(V^X)}{\sqrt{-p'(v_-)}}(\sqrt{-p'(v_-)}\phi+\psi)dy.
\end{align*}
The last three terms on the right hand side can be controlled as follows:
\begin{align*}
  &\delta\int_0^1\left(sa^X-2\sqrt{-p'(v_-)}-\frac{a^Xp'(V^X)}{\sqrt{-p'(v_-)}}\right)\omega dy\leq C\delta(\delta+\lambda)\|\omega\|_{L^2(0,1)}, \\
  &\delta\int_0^1a^X\frac{p'(V^X)}{\sqrt{-p'(v_-)}}(\sqrt{-p'(v_-)}\phi+\psi)dy\leq C\sqrt{\frac{\delta}{k}}\left(\int_{\bbR}a^X_{\xi}(\sqrt{-p'(v_-)}\phi+\psi)^2dx\right)^{\frac12}.
\end{align*}
According to the inequality $(A+B)^2\geq \frac12A^2-B^2$, we obtain
\begin{align*}
  \left(\int_{\bbR}a^X(U^X_{\xi}\psi-p'(V^X)V^X_{\xi}\phi)dx\right)^2\geq&\, -2\delta^2p'(v_-)\left(\int_0^1\omega dy\right)^2
  -C\delta^2(\delta+\lambda)^2\int_0^1\omega^2dy\\
  &\,-C\frac{\delta}{k}\int_{\bbR}a^X_{\xi}(\sqrt{-p'(v_-)}\phi+\psi)^2dx.
\end{align*}
Hence,
\begin{align}\label{Y-1}
  &\,\frac{1}{2}\int_{\mathfrak{G}}\dot{X}(t)\int_{\bbR}a^X(U^X_{\xi}\psi-p'(V^X)V^X_{\xi}\phi)dxdt=-\frac{\delta(-p'(v_-))^{\frac32}}{2p''(v_-)}\int_{\mathfrak{G}}|\dot{X}(t)|^2dt\notag\\
  =&\,-\frac{p''(v_-)}{2\delta(-p'(v_-))^{\frac32}}\int_{\mathfrak{G}}\left(\int_{\bbR}a^X(U^X_{\xi}\psi-p'(V^X)V^X_{\xi}\phi)dx\right)^2dt\notag\\
  \leq&\, -\delta\frac{p''(v_-)}{\sqrt{-p'(v_-)}}\int_{\mathfrak{G}}\left(\int_0^1\omega dy\right)^2dt+C\delta(\delta+\lambda)^2\int_{\mathfrak{G}}\int_0^1\omega^2dydt+\frac{C}{k}\int_{\mathfrak{G}}\int_{\bbR}a^X_{\xi}(\sqrt{-p'(v_-)}\phi+\psi)^2dxdt.
\end{align}

The second term in $Y(v,u)$ can be estimated as follows:
\begin{align}\label{Y-2}
  \int_{\mathfrak{G}}\dot{X}(t)\int_{\bbR}a^X_{\xi}\eta((v,u)|(V^X,U^X))dxdt\leq&\, Ck\int_{\mathfrak{G}}|\dot{X}(t)|\int_{\bbR}V^X_{\xi}(\psi^2+\phi^2)dxdt\notag\\
  \leq&\, \frac{\delta(-p'(v_-))^{\frac32}}{4p''(v_-)}\int_{\mathfrak{G}}|\dot{X}(t)|^2dt+C\frac{k^2}{\delta}\int_{\mathfrak{G}}\left(\int_{\bbR}V^X_{\xi}(\psi^2+\phi^2)dx\right)^2dt\notag\\
  \leq&\, \frac{\delta(-p'(v_-))^{\frac32}}{4p''(v_-)}\int_{\mathfrak{G}}|\dot{X}(t)|^2dt+C\eta^2 k^2\int_{\mathfrak{G}}\int_{\bbR}V^X_{\xi}(\psi^2+\phi^2)dxdt,
\end{align}
where we have used the fact that
\begin{equation*}
  \sup_{t\in\mathfrak{G}}\int_{\bbR}V^X_{\xi}(\psi^2+\phi^2)dx\leq \sup_{t\in\mathfrak{G}}\|(\phi,\psi)\|_{L^{\infty}}^2\int_{\bbR}V^X_{\xi}dx\leq \delta\eta^2.
\end{equation*}

Taking $k$ suitably large, for small enough $\eta$ and $\delta$, we derive from \eqref{J-bad-1}, \eqref{J-good-1} and \eqref{Y-1} that
\begin{align}\label{good-term}
  &\,s\int_{\mathfrak{G}}\int_{\bbR}a^XV^X_{\xi}p(v|V^X)dxdt-\frac34\int_{\mathfrak{G}}\int_{\bbR}\frac{\mu}{v}a^X\psi_x^2dxdt\\
  &\,\quad-\frac{\sqrt{-p'(v_-)}}{2}\int_{\mathfrak{G}}\int_{\bbR}a^X_{\xi}(\sqrt{-p'(v_-)}\phi+\psi)^2dxdt-\frac{\delta(-p'(v_-))^{\frac32}}{2p''(v_-)}\int_{\mathfrak{G}}|\dot{X}(t)|^2dt\notag\\
  \leq&\, -C\left(\delta\int_{\mathfrak{G}}\int_0^1\omega^2dydt+\delta\int_{\mathfrak{G}}\left(\int_0^1\omega\right)^2dt+\int_{\mathfrak{G}}\int_{\bbR}a^X_{\xi}(\sqrt{-p'(v_-)}\phi+\psi)^2dxdt\right)\notag\\
  &\,\quad +C(\delta+\eta)\int_{\mathfrak{G}}\int_{\bbR}V^X_{\xi}(\phi^2+\psi^2)dxdt\notag\\
  \leq&\, -c_1\int_{\mathfrak{G}}\int_{\bbR}V^X_{\xi}(\phi^2+\psi^2)dxdt,
\end{align}
where we have used the fact that
\begin{equation*}
  \int_{\mathfrak{G}}\int_{\bbR}V^X_{\xi}(\phi^2+\psi^2)dxdt\leq \delta\int_{\mathfrak{G}}\int_0^1\omega^2dydt+\frac{1}{k}\int_{\mathfrak{G}}\int_{\bbR}a^X_{\xi}(\sqrt{-p'(v_-)}\phi+\psi)^2dxdt.
\end{equation*}

Combining \eqref{J-bad-2-4}-\eqref{J-bad-7}, \eqref{Y-2}, \eqref{good-term} with \eqref{rela-entropy-integ} yields
\begin{align}\label{good-est-final}
  &\,\int_{\bbR}a^X\eta((v,u)|(V^X,U^X))dx+\frac{1}{8}\int_{\mathfrak{G}}\int_{\bbR}\frac{\mu}{v}a^X\psi_x^2dxdt+\int_{\mathfrak{B}}\int_{\bbR}\frac{\mu}{v}a^X\psi_x^2dxdt\notag\\
  &\,\quad +c_1\int_{\mathfrak{G}}\int_{\bbR}V^X_{\xi}(\phi^2+\psi^2)dxdt+\frac{\delta(-p'(v_-))^{\frac32}}{4p''(v_-)}\int_{\mathfrak{G}}|\dot{X}(t)|^2dt+\frac{(-p'(v_-))^{\frac32}}{p''(v_-)}\delta\int_{\mathfrak{B}}|\dot{X}(t)|^2dt\notag\\
  \leq&\, \int_{\bbR}a^X\eta((v,u)|(V^X,U^X))(0,x)dx+C(k(\delta+\eta)+(1+k^2)\delta^2+k^2\eta^2)\int_{\mathfrak{G}}\int_{\bbR}V^X_{\xi}(\phi^2+\psi^2)dxdt\notag\\
  &\,-\int_{\mathfrak{B}}\int_{\bbR}\mu a^XU^X_{\xi}\psi_x\left(\frac{1}{v}-\frac{1}{V^X}\right)dxdt-\int_{\mathfrak{B}}\int_{\bbR}\mu a^X_{\xi}\psi\left(\frac{\psi_x}{v}+U^X_{\xi}\left(\frac{1}{v}-\frac{1}{V^X}\right)\right)dxdt\notag\\
  &\,+\int_{\mathfrak{B}}\int_{\bbR}a^X_{\xi}\left[\psi(p(v)-p(V^X))-s\eta((v,u)|(V^X,U^X))\right]dxdt+s\int_{\mathfrak{B}}\int_{\bbR}a^XV^X_{\xi}p(v|V^X)dxdt\notag\\
  &\,-\int_{\mathfrak{B}}\dot{X}(t)\int_{\bbR}a^X_{\xi}\eta((v,u)|(V^X,U^X))dxdt.
\end{align}
Taking $k(\delta+\eta)\ll1$, we can absorb the first term on the right hand side of \eqref{good-est-final}.
\medbreak

\noindent{\bf 2. Bad time estimates}

Next, we derive rough bounds for the rest terms on the right hand side of \eqref{good-est-final}.
We define the terms considering in bad time set by $I_{j}$ $(j=1\sim 7)$. Notice that
$I_5$ have a good sign which can be omitted.
Using Cauchy inequality, we obtain
\begin{align*}
  |I_1|\leq&\, \frac{1}{4}\int_{\mathfrak{B}}\int_{\bbR}\frac{\mu}{v} a^X\psi_x^2dxdt+C\int_{\mathfrak{B}}\int_{\bbR}\mu va^x(U^X_{\xi})^2\left(\frac{1}{v}-\frac{1}{V^X}\right)^2dxdt\\
  \leq&\, \frac{1}{4}\int_{\mathfrak{B}}\int_{\bbR}\frac{\mu}{v} a^X\psi_x^2dxdt+C(M_1)\delta^2\int_{\mathfrak{B}}\int_{\bbR}V^X_{\xi}(\phi^2+\psi^2)dxdt,\\
  |I_2|\leq&\, \frac{1}{4}\int_{\mathfrak{B}}\int_{\bbR}\frac{\mu}{v} a^X\psi_x^2dxdt+C\int_{\mathfrak{B}}\int_{\bbR}\mu \frac{(a_{\xi}^X)^2}{a^Xv}\psi^2dxdt\\
  \leq&\, \frac{1}{4}\int_{\mathfrak{B}}\int_{\bbR}\frac{\mu}{v} a^X\psi_x^2dxdt+C(M_1)k^2\delta^2\int_{\mathfrak{B}}\int_{\bbR}V^X_{\xi}(\phi^2+\psi^2)dxdt,\\
  |I_3|\leq&\, C(M_1)k\int_{\mathfrak{B}}\int_{\bbR}(V^X_{\xi})^2|\phi\psi|dxdt\leq C(M_1)k\delta^2\int_{\mathfrak{B}}\int_{\bbR}V^X_{\xi}(\phi^2+\psi^2)dxdt,\\
  |I_4|\leq&\, C(M_1)k\int_{\mathfrak{B}}\int_{\bbR}V^X_{\xi}|\phi\psi|dxdt\leq C(M_1)k\int_{\mathfrak{B}}\int_{\bbR}V^X_{\xi}(\phi^2+\psi^2)dxdt,\\
  |I_6|\leq&\, C(M_1)\int_{\mathfrak{B}}\int_{\bbR}V^X_{\xi}|\phi|^2dxdt\leq C(M_1)\int_{\mathfrak{B}}\int_{\bbR}V^X_{\xi}(\phi^2+\psi^2)dxdt,\\
  |I_7|\leq&\, \frac{\delta(-p'(v_-))^{\frac32}}{2p''(v_-)}\int_{\mathfrak{B}}|\dot{X}(t)|^2dt+C\delta^{-1}\int_{\mathfrak{B}}\left(\int_{\bbR}a^X_{\xi}\eta((v,u)|(V^X,U^X))dx\right)^2dt\\
  \leq&\, \frac{\delta(-p'(v_-))^{\frac32}}{2p''(v_-)}\int_{\mathfrak{B}}|\dot{X}(t)|^2dt+C(M_1)\frac{k^2}{\delta}\int_{\mathfrak{B}}\left(\int_{\bbR}V^X_{\xi}(\phi^2+\psi^2)dx\right)^2dt\\
  \leq&\, \frac{\delta(-p'(v_-))^{\frac32}}{2p''(v_-)}\int_{\mathfrak{B}}|\dot{X}(t)|^2dt+C(M_1)k^2\delta\int_{\mathfrak{B}}\int_{\bbR}a^X\eta((v,u)|(V^X,U^X))dx\int_{\bbR}V^X_{\xi}(\phi^2+\psi^2)dxdt\\
  \leq&\, \frac{\delta(-p'(v_-))^{\frac32}}{2p''(v_-)}\int_{\mathfrak{B}}|\dot{X}(t)|^2dt+C(M_1)k^2\int_{\mathfrak{B}}\int_{\bbR}V^X_{\xi}(\phi^2+\psi^2)dxdt.
\end{align*}
Inserting the above calculations into \eqref{good-est-final} yields the following estimates
\begin{align*}
  &\,\int_{\bbR}a^X\eta((v,u)|(V^X,U^X))dx+\frac{1}{8}\int_{\mathfrak{G}}\int_{\bbR}\frac{\mu}{v}a^X\psi_x^2dxdt+\frac12\int_{\mathfrak{B}}\int_{\bbR}\frac{\mu}{v}a^X\psi_x^2dxdt\notag\\
  &\,\quad +c_2\int_{\mathfrak{G}}\int_{\bbR}V^X_{\xi}(\phi^2+\psi^2)dxdt+\frac{\delta(-p'(v_-))^{\frac32}}{4p''(v_-)}\int_{\mathfrak{G}}|\dot{X}(t)|^2dt+\frac{(-p'(v_-))^{\frac32}}{2p''(v_-)}\delta\int_{\mathfrak{B}}|\dot{X}(t)|^2dt\notag\\
  \leq&\,\int_{\bbR}a^X\eta((v,u)|(V^X,U^X))(0,x)dx+C(M_1)(1+k^2)\int_{\mathfrak{B}}\int_{\bbR}V^X_{\xi}(\phi^2+\psi^2)dxdt
\end{align*}
Combining the same terms with respect to good time set and bad time set, we obtain
\begin{align}\label{gb-est-final}
  &\,\int_{\bbR}a^X\eta((v,u)|(V^X,U^X))dx+\frac{1}{8}\int_0^T\int_{\bbR}\frac{\mu}{v}a^X\psi_x^2dxdt\notag\\
  &\,\quad +c_2\int_0^T\int_{\bbR}V^X_{\xi}(\phi^2+\psi^2)dxdt+\frac{\delta(-p'(v_-))^{\frac32}}{4p''(v_-)}\int_0^T|\dot{X}(t)|^2dt\notag\\
  \leq&\,\int_{\bbR}a^X\eta((v,u)|(V^X,U^X))(0,x)dx+(C(M_1)(1+k^2)+c_2)\int_0^T1_{\mathfrak{B}}\int_{\bbR}V^X_{\xi}(\phi^2+\psi^2)dxdt.
\end{align}
Notice that
\begin{align*}
  &1_{\mathfrak{B}(t)}\leq \eta^{-4}\|(\phi,\psi)\|_{L^{\infty}}^4\leq C(M_1)\eta^{-4}N_2\|(\phi_x,\psi_x)\|_{L^2}^2,\\
  &\int_{\bbR}V^X_{\xi}(\phi^2+\psi^2)dx\leq C(M_1)\delta^2 N_2.
\end{align*}
We finally have
\begin{align*}
  &\,\int_{\bbR}a^X\eta((v,u)|(V^X,U^X))dx+c_0\int_0^t\left(\int_{\bbR}\frac{\mu}{v}a^X\psi_x^2dx+\int_{\bbR}V^X_{\xi}(\phi^2+\psi^2)dx+\delta|\dot{X}(t)|^2\right)d\tau\\
  \leq&\, \int_{\bbR}a^X\eta((v,u)|(V^X,U^X))(0,x)dx+C(M_1)(1+k^2)\eta^{-4}\delta^2N^2_2\int_0^t\|(\phi_x,\psi_x)\|_{L^2}^2d\tau.
\end{align*}
\end{proof}

The next lemma is focused on the estimates of $\int_0^t\|(\phi_x,\psi_x)\|_{L^2}^2d\tau$.
\begin{lemma}\label{lem-first-order}
Under the assumptions in Proposition \ref{proposition-o-f}, the following estimate holds
\begin{align}\label{first-order-est}
  &\,\|(\log v-\log V^X)_x\|_{L^2}^2+\int_0^T\|(\phi_x,\psi_x)\|_{L^2}^2dt\notag\\
  \leq&\, C(M_1)\left(A_0+\int_{\bbR}a^X\eta((v,u)|(V^X,U^X))dx+\int_0^T\int_{\bbR}\frac{\mu}{v}a^X\psi_x^2dxdt+\delta^2\int_0^T\int_{\bbR}V^X_{\xi}(\phi^2+\psi^2)dxdt\right),
\end{align}
where $A_0=1+\int_{\bbR}a^X\eta((v,u)|(V^X,U^X))(0,x)dx+\|(\log v-\log V^X)_x(0,x)\|_{L^2}^2$.
\end{lemma}
\begin{proof}
We consider the following functional
\begin{equation*}
  \mathcal{B}(t)=\frac{\mu^2}{2}\|(\log v-\log V^X)_x\|_{L^2}^2-\int_{\bbR}\mu \psi(\log v-\log V^X)_xdx,
\end{equation*}
which satisfies
\begin{align}\label{B-lb}
  \mathcal{B}(t)\geq&\, \frac{\mu^2}{4}\|(\log v-\log V^X)_x\|_{L^2}^2-\|\psi\|_{L^2}^2,\\
  \mathcal{B}(0)\leq&\, \frac{3\mu^2}{4}\|(\log v-\log V^X)_x(0,x)\|_{L^2}^2+\|\psi(0)\|_{L^2}^2.\label{B0-ub}
\end{align}
Direct calculations show that
\begin{align}\label{B}
  \frac{d}{dt}\mathcal{B}(t)=&\,\int_{\bbR}\mu(\log v-\log V^X)_x(p(v)-p(V^X))_xdx-\int_{\bbR}\mu\psi\left(\frac{u_x}{v}-\frac{U^X_{\xi}}{V^X}\right)_xdx\notag\\
  &\,+\dot{X}(t)\left[\mu^2\int_{\bbR}(\log v-\log V^X)_x\left(\frac{V^X_{\xi}}{V^X}\right)_{\xi}dx-\mu\int_{\bbR}(\log v-\log V^X)_xU^X_{\xi}dx-\mu\int_{\bbR}\psi\left(\frac{V^X_{\xi}}{V^X}\right)_{\xi}dx\right]\notag\\
  :=&\,\sum_{i=1}^5\mathcal{B}_i(t).
\end{align}

Notice that
\begin{equation*}
  (p(v)-p(V^X))_x=-\gamma p(v)(\log v-\log V^X)_x-\gamma(p(v)-p(V^X))\frac{V^X_{\xi}}{V^X},
\end{equation*}
The first term on the right hand side of \eqref{B} can be controlled as follows:
\begin{align*}
  \mathcal{B}_1=&\,- \gamma\mu \int_{\bbR}p(v)(\log v-\log V^X)_x^2dx-\gamma\mu\int_{\bbR}(\log v-\log V^X)_x(p(v)-p(V^X))\frac{V^X_{\xi}}{V^X}dx\\
  \leq&\,-\frac{\gamma\mu}{2}\int_{\bbR}p(v)(\log v-\log V^X)_x^2dx+C(M_1)\delta^2\int_{\bbR}V^X_{\xi}\phi^2dx.
\end{align*}
For the term $\mathcal{B}_2$, we apply Cauchy's inequality to get
\begin{equation*}
  \mathcal{B}_2\leq C(M_1)\left(\int_{\bbR}\mu\frac{\psi_x^2}{v}dx+\int_{\bbR}(U_{\xi}^X)^2\phi^2dx\right)\leq C(M_1)\left(\int_{\bbR}\mu\frac{\psi_x^2}{v}dx+\delta^2\int_{\bbR}V_{\xi}^X\phi^2dx\right).
\end{equation*}
For the last three terms, we can obtain similarly
\begin{align*}
  |\mathcal{B}_3|+|\mathcal{B}_4|\leq&\, \frac{\gamma\mu}{4}\int_{\bbR}p(v)(\log v-\log V^X)_x^2dx+C(M_1)|\dot{X}(t)|^2\delta^3,\\
  \mathcal{B}_5=&\,\dot{X}(t)\int_{\bbR}\mu\psi_x\frac{V^X_{\xi}}{V^X}dx\leq \int_{\bbR}\mu\frac{\psi_x^2}{v}dx+C(M_1)|\dot{X}(t)|^2\delta^3.
\end{align*}

Inserting the above estimates into \eqref{B}, combining with \eqref{B-lb} and \eqref{B0-ub}, we get
\begin{align}\label{phi-x}
  &\,\|(\log v-\log V^X)_x\|_{L^2}^2+\frac{\gamma\mu}{4}\int_0^T\int_{\bbR}p(v)(\log v-\log V^X)_x^2dxdt\notag\\
  \leq&\, C(M_1)\left(A_0+\int_{\bbR}a^X\eta((v,u)|(V^X,U^X))dx+\int_0^T\int_{\bbR}\frac{\mu}{v}a^X\psi_x^2dxdt+\delta^2\int_0^T\int_{\bbR}V^X_{\xi}(\phi^2+\psi^2)dxdt\right),
\end{align}
Since $\phi_x$ satisfies
\begin{equation*}
\phi_x=v(\log v-\log V^X)_x+\frac{V^X_{\xi}}{V^X}\phi,
\end{equation*}
we have
\begin{equation*}
  \|\phi_x\|_{L^2}^2\leq C(M_1)\int_{\bbR}p(v)(\log v-\log V^X)_x^2dx+C\delta^2\int_{\bbR}V^X_{\xi}\phi^2dx.
\end{equation*}
Inserting above estimate into \eqref{phi-x}, we obtain the desired estimate \eqref{first-order-est}.
\end{proof}

Combining Lemma \ref{lem-gb-time} with \eqref{first-order-est}, we have
\begin{lemma}\label{lem-L2-close}
Under the assumptions in Proposition \ref{proposition-o-f}. Suppose that
\begin{equation}\label{assum-l2}
    N_2=4A_0,\quad C(M_1)(1+k^2)\eta^{-4}\delta^2N^2_2\leq \min\left\{\frac14,\frac{c_0}{4}\right\},
\end{equation}
then the following estimates hold
\begin{align}\label{L2-est-1}
  &\int_{\bbR}a^X\eta((v,u)|(V^X,U^X))dx\leq 2A_0,\\
  &\int_{\bbR}a^X\eta((v,u)|(V^X,U^X))dx+\int_0^T\left(\int_{\bbR}\frac{\mu}{v}a^X\psi_x^2dx+\int_{\bbR}V^X_{\xi}(\phi^2+\psi^2)dx+\delta|\dot{X}(t)|^2\right)dt
  \leq C.\label{L2-est-2}
\end{align}
\end{lemma}

\subsection{Uniformly bounds for specific volume}

\begin{lemma}\label{lem-v-bd}
Under the assumptions in Proposition \ref{proposition-o-f}, there exist two positive constants $\underline{v}$ and $\overline{v}$ which only depend on $A_0$, $N_1^0$ and $M_1^0$ such that
\begin{equation}\label{v-bd}
  \underline{v}\leq v(t,x)\leq \overline{v}.
\end{equation}
\end{lemma}
\begin{proof}
We divide the proof into two steps: In the first step, we will derive that the mean value of the specific volume is uniformly bounded in any interval $I$ of length $1$.
In the second step, we will compare $v(t,x)$ with $\int_{I}v(t,x)dx$ in any interval $I$ of length $1$ containing $x$.

\smallskip
{\it Step 1.} Since $V^X\in[\frac{v_+}{2}, v_+]$, there exists generically constant $C$ such that
\begin{equation*}
  0<Q(v|v_+)\leq CQ(v|V^X)+C.
\end{equation*}
Then, applying the Jensen's inequality, we obtain
\begin{equation*}
  Q\left(\int_{I}vdx|v_+\right)\leq \int_{I}Q(v|v_+)dx\leq C(A_0+1).
\end{equation*}
Since $\lim_{v\rightarrow0}Q(\int_{I}vdx|v_+)=\lim_{v\rightarrow\infty}Q(\int_{I}vdx|v_+)=\infty$, there exist two constants $c(A_0)$ and $C(A_0)$ such that
\begin{equation}\label{mean-value-bd}
  0<c(A_0)\leq \int_{I}v(t,x)dx\leq C(A_0).
\end{equation}
Hence, there exists $y_{I}\in I$ such that $v(t,y_I)=\int_{I}v(t,x)dx$.

\smallskip
{\it Step 2.} For a interval $I$ of length $1$ that contains any $x\in\bbR$, we introduce the following quantity:
\begin{equation}\label{L}
  L(t,x;y_I)=\mu\log v(t,x)-\mu\log v(t,y_I)-\int_{y_I}^x u(t,z)dz, \quad {\rm for\ any\ } |x-y_I|\leq 1,
\end{equation}
which satisfies
\begin{equation*}
  L(t,x;y_I)_t=p(v(t,x))-p(v(t,y_I)),
\end{equation*}
and
\begin{align*}
  |L(0,x;y_I)|\leq \left|\mu\int_{y_I}^x(\log v(0,x))_xdx\right|+\left|\int_{y_I}^x u(0,z)dz\right|
  \leq (L_1(M_1^0,N_1^0)+L_2(A_0))\sqrt{|x-y_I|}.
\end{align*}
Applying H\"{o}lder's inequality, we have
\begin{equation*}
  \left|\int_{y_I}^x u(t,z)dz\right|\leq \sqrt{8A_0}\sqrt{|x-y_I|}+\|U^X\|_{L^{\infty}}|x-y_I|\leq L_2(A_0)\sqrt{|x-y_I|}.
\end{equation*}

If $L(t,x;y_I)>(L_1(M_1^0,N_1^0)+L_2(A_0))\sqrt{x-y_I}$, according to \eqref{L}, we get
\begin{equation*}
  \mu\log\left(\frac{v(t,x)}{v(t,y_I)}\right)>L(t,x;y_I)-L_2(A_0)\sqrt{|x-y_I|}>L_1(M_1^0,N_1^0)\sqrt{|x-y_I|}>0,
\end{equation*}
which yields $v(t,x)>v(t,y_I)$, hence $L(t,x;y_I)_t<0$. We then have
\begin{equation*}
  \mu\log\left(\frac{v(t,x)}{v(t,y_I)}\right)\leq L(0,x;y_I)+L_2(A_0)\sqrt{|x-y_I|}\leq L_1(M_1^0,N_1^0)+2L_2(A_0).
\end{equation*}

If $L(t,x;y_I)<-(L_1(M_1^0,N_1^0)+L_2(A_0))\sqrt{|x-y_I|}$, similarly, we have
\begin{equation*}
\mu\log\left(\frac{v(t,x)}{v(t,y_I)}\right)<-L_1(M_1^0,N_1^0)\sqrt{|x-y_I|}<0
\end{equation*}
which means $L(t,x;y_I)_t>0$. Hence, we obtain
\begin{equation*}
  \mu\log\left(\frac{v(t,x)}{v(t,y_I)}\right)\geq L(0,x;y_I)-L_2(A_0)\sqrt{|x-y_I|}\geq -(L_1(M_1^0,N_1^0)+2L_2(A_0)).
\end{equation*}

If $|L(t,x;y_I)|\leq (L_1(M_1^0,N_1^0)+L_2(A_0))\sqrt{|x-y_I|}$, then
\begin{equation*}
  \left|\mu\log\left(\frac{v(t,x)}{v(t,y_I)}\right)\right|\leq (L_1(M_1^0,N_1^0)+2L_2(A_0))\sqrt{|x-y_I|}.
\end{equation*}

The above three cases yields the final upper and lower bounds on the specific volume $v(t,x)$. This completes the proof of this lemma.
\end{proof}

\subsection{Higher order estimates}
\begin{lemma}\label{lem-high-order}
Under the assumptions in Proposition \ref{proposition-o-f}, we have
\begin{align}\label{high-order-est}
  \sup_{t\in[0,T]}\|(\phi,\psi)\|_{H^1}^2+\int_0^T(\|\phi_x\|_{L^2}^2+\|\psi_x\|_{H^1}^2)dt\leq C.
\end{align}
\end{lemma}
\begin{proof}
Firstly, we focus on $\phi_x(t,x)$.
Applying Lemma \ref{lem-first-order}, \eqref{L2-est-2} and Lemma \ref{lem-v-bd}, we obtain
\begin{equation}\label{phi-1st-est}
  \sup_{t\in[0,T]}(\|\phi\|_{H^1}^2+\|\psi\|^2_{L^2})+\int_0^T(\|\phi_x\|_{L^2}^2+\|\psi_x\|_{L^2}^2)dt\leq C.
\end{equation}

Recall that
\begin{equation*}
  \psi_t+(p(v)-p(V^X))_x=\left(\mu\frac{\psi_x}{v}-\mu U^X_{\xi}\frac{\phi}{vV^X}\right)_x+\dot{X}(t)U^X_{\xi}.
\end{equation*}
Multiplying both sides by $-\psi_{xx}$ and integrating over $\bbR$ with respect to $x$, we have
\begin{align*}
  \frac{1}{2}\frac{d}{dt}\|\psi_x\|_{L^2}^2+\int_{\bbR}\frac{\mu}{v}\psi_{xx}^2dx=&\,\int_{\bbR}(p(v)-p(V^X))_x\psi_{xx}dx
  +\int_{\bbR}\mu\frac{v_x}{v^2}\psi_x\psi_{xx}dx\\
  &\,+\int_{\bbR}\left(\mu U^X_{\xi}\frac{\phi}{vV^X}\right)_x\psi_{xx}dx-\int_{\bbR}\dot{X}(t)U^X_{\xi}\psi_{xx}dx.
\end{align*}
Applying Lemma \ref{lem-v-bd}, we get
\begin{align*}
  \|(p(v)-p(V^X))_x\|_{L^2}^2\leq&\, C\left(\|\phi_x\|_{L^2}^2+\delta^2\int_{\bbR}V_{\xi}^X\phi^2dx\right),\\
  \left\|\left(\mu U^X_{\xi}\frac{\phi}{vV^X}\right)_x\right\|_{L^2}^2\leq&\, C\left(\delta^4\|\phi_x\|_{L^2}^2+\delta^2\int_{\bbR}V_{\xi}^X\phi^2dx\right),\\
  \|\dot{X}(t)U^X_{\xi}\|_{L^2}^2\leq&\, C\delta^3|\dot{X}(t)|^2.
\end{align*}
Then, \eqref{phi-1st-est} and H\"{o}lder's inequality yield
\begin{align*}
  \int_{\bbR}\mu\frac{v_x}{v^2}\psi_x\psi_{xx}dx\leq&\, C\|v_x\|_{L^2}\|\psi_x\|_{L^{\infty}}\|\psi_{xx}\|_{L^2}\\
  \leq&\, \epsilon\|\psi_{xx}\|_{L^2}^2+C(\epsilon)(\|\phi_x\|_{L^2}^4+\|V_{\xi}^X\|_{L^2}^4)\|\psi_x\|_{L^{2}}^2\\
  \leq&\, \epsilon\|\psi_{xx}\|_{L^2}^2+C(\epsilon)\|\psi_x\|_{L^{2}}^2.
\end{align*}

Finally, with \eqref{L2-est-2}, \eqref{phi-1st-est} and H\"{o}lder's inequality, we can derive
\begin{align*}
  \|\psi_x\|_{L^2}^2+\int_0^T\int_{\bbR}\frac{\mu}{v}\psi_{xx}^2dx\leq \|\psi_x(0,x)\|_{L^2}^2+C\int_0^T\left(\|\phi_x\|_{L^2}^2+\delta^2\int_{\bbR}V_{\xi}^X\phi^2dx+\delta^3|\dot{X}(t)|^2+\|\psi_x\|_{L^{2}}^2\right)dt\leq C.
\end{align*}
This combines with \eqref{phi-1st-est} completes the proof of this lemma.
\end{proof}

\section{The proof of Theorem \ref{theorem1}}
In this section, we prove Theorem \ref{theorem1}. The global existence can be derived by the local existence in Proposition \ref{local-exist}, and
the a priori estimates in Lemmas \ref{lem-gb-time}-\ref{lem-high-order} combining with the continuity argument below. For brevity, we set
\begin{equation*}
  E(t):=\int_{\bbR}a^X\eta((v,u)|(V^X,U^X))dx\leq N_2,\quad {\rm with\ } N_2=4A_0.
\end{equation*}
First, we choose a weak shock wave such that $0<\delta\leq\delta_1\leq \min\{1,\frac{v_+}{2}\}$, then the viscous shock profiles are uniformly bounded.
Next, we choose $k$ sufficiently large, then assume that the parameter $\eta$ of good time is sufficiently small such that $k\eta\ll1$, and further restrict
$\delta_1$ such that $k\delta_1\ll1$  as in Lemma \ref{lem-gb-time}.

Then, under the assumption $E(t)\leq N_2$. Lemma \ref{lem-v-bd} gives $\underline{v}\leq v\leq \overline{v}$. Hence, we set
\begin{equation*}
  M_{\ast}=\max\left\{1, M_1^0,\underline{v}^{-1}, \overline{v}, v_+, \frac{2}{v_+}\right\},\quad {\rm and\ } \quad M_1=4M_{\ast}.
\end{equation*}

With above preparation, under \eqref{L2-est-2} and Lemma \ref{lem-v-bd}, we apply Lemma \ref{lem-high-order} to get
\begin{align*}
  \sup_{t\in[0,T]}\|(\phi,\psi)\|_{H^1}^2+\int_0^T(\|\phi_x\|_{L^2}^2+\|\psi_x\|_{H^1}^2)dt\leq C_{\ast}.
\end{align*}
The constant $C_{\ast}$ is independent of $N_1$. Then, choose $N_1=\max\{1, N_1^0, \sqrt{C_{\ast}}\}$.

Finally, choose the strength of shock wave $0<\delta_0\leq \delta_1$ sufficiently small such that the smallness assumptions in Lemma \ref{lem-gb-time} are
satisfied for $0<\delta\leq\delta_0$ and that
\begin{equation*}
  \delta N_2\leq \epsilon_0,\quad C(M_1)(1+k^2)\eta^{-4}\delta^2N^2_2\leq \min\left\{\frac14,\frac{c_0}{4}\right\}.
\end{equation*}
Here the positive constant $\epsilon_0$ is sufficiently small and the choice of $\delta$ depends only on the initial data and the fixed parameters.

Let $T_{max}$ be the maximal existence time of the local strong solution together with its shift. By Proposition \ref{local-exist}, the following bootstrap conditions hold
\begin{equation*}
  M_1^{-1}\leq v\leq M_1,\quad \|(\phi,\psi)\|_{H^1}\leq N_1,\quad E(t)\leq N_2.
\end{equation*}
Let $T_{\ast}$ be the supremum of all $T<T_{max}$ for which these conditions hold on $[0, T]$. For every $T<T_{\ast}$,
Lemmas \ref{lem-L2-close}-\ref{lem-high-order} yield the strict improvements
\begin{equation*}
  \frac{4}{M_1}\leq v\leq \frac{M_1}{4},\quad \|(\phi,\psi)\|_{H^1}\leq \frac{N_1}{4},\quad E(t)\leq \frac{N_2}{2}.
\end{equation*}
Hence $T_{\ast}=T_{max}$, and all the preceding a priori estimates hold uniformly on every compact subinterval of $[0, T_{max})$.
Set $X_{\ast}=\lim_{t\rightarrow T_{max}} X(t)$, and
\begin{equation*}
  (v_{\ast},u{\ast})(x)=(\phi,\psi)(x)+(V,U)(x-sT_{max}-X_{\ast}).
\end{equation*}
Applying Proposition \ref{local-exist} at $t=T_{max}$, we can extend the local strong solution beyond $T_{max}$ which contradicts maximality. Therefore, $T_{max}=\infty$.

The large time behavior is established below.
Since $\phi_{tx}=\psi_{xx}+\dot{X}(t)V^X_{\xi\xi}$, we apply H\"{o}lder's inequality to get
\begin{align*}
  \left|\frac{d}{dt}\|\phi_x\|_{L^2}^2\right|\leq&\, 2\|\phi_x\|_{L^2}(\|\psi_{xx}\|_{L^2}+|\dot{X}(t)|\|V^X_{\xi\xi}\|_{L^2})\\
  \leq&\, 2\|\phi_x\|_{L^2}(\|\psi_{xx}\|_{L^2}+\delta^{\frac52}|\dot{X}(t)|),\\
  \left|\frac{d}{dt}\|\psi_x\|_{L^2}^2\right|\leq&\, 2\|\psi_{xx}\|_{L^2}\|\psi_t\|_{L^2}.
\end{align*}
Hence, combining with Lemma \ref{lem-high-order}, we obtain
\begin{align*}
  \int_0^{\infty}\left(\|\phi_x\|_{L^2}^2+\left|\frac{d}{dt}\|\phi_x\|_{L^2}^2\right|\right)dt
  +\int_0^{\infty}\left(\|\psi_x\|_{L^2}^2+\left|\frac{d}{dt}\|\psi_x\|_{L^2}^2\right|\right)dt\leq C,
\end{align*}
which yields
\begin{equation*}
  \lim_{t\rightarrow\infty}\|(\phi,\psi)_x\|_{L^{2}}=0,\quad \lim_{t\rightarrow\infty}\|(\phi,\psi)\|_{L^{\infty}}=0.
\end{equation*}
We finish the proof of Theorem \ref{theorem1}.

\medbreak


\noindent{\bf AI Declaration:}\, The authors acknowledge the use of AI tools for exploratory discussions and literature searches.
The authors independently worked and verified all theorem statements, proofs, calculations and take full responsibility for the mathematical content.

\bibliographystyle{siam}


 \end{document}